\documentclass[12pt]{article}

\usepackage{amsmath}
\usepackage{amsthm}
\usepackage{amssymb}
\usepackage{latexsym}
\usepackage{pstricks}
\usepackage{graphicx, pst-plot, pst-node, pst-text, pst-tree}

\theoremstyle{plain}
\newtheorem{thm}{Theorem}[section]

\newtheorem{lem}[thm]{Lemma}
\newtheorem{prop}[thm]{Proposition}

\newtheorem{question}[thm]{Question}
\theoremstyle{definition}

\newtheorem{exa}[thm]{Example}
\theoremstyle{remark}

\newcounter{todocounter}

\makeatletter
\newfont{\footsc}{cmcsc10 at 8truept}
\newfont{\footbf}{cmbx10 at 8truept}
\newfont{\footrm}{cmr10 at 10truept}
\title{Wreath Products of Permutation Classes}
\author{Robert Brignall\\[-5pt]
\small School of Mathematics and Statistics\\[-5pt]
\small University of St~Andrews\\[-5pt]
\small St~Andrews, Fife, Scotland\\[-5pt]
\small \texttt{robertb@mcs.st-and.ac.uk}\\[-5pt]
\small \texttt{http://turnbull.mcs.st-and.ac.uk/\~{}robertb}}

\date{\today \\[6pt]
        \begin{flushleft}
        \small Key Words: permutation class, restricted permutation, basis, wreath product\\[6pt]
        \end{flushleft}
}

\long\def\symbolfootnote[#1]#2{\begingroup%
\def\thefootnote{\fnsymbol{footnote}}\footnote[#1]{#2}\endgroup}

\begin{document}
\maketitle

\newcommand{\OEISlink}[1]{{#1}}
\newcommand{\OEISref}{{OEIS}~\cite{OEIS}}
\newcommand{\OEIS}[1]{(Sequence \OEISlink{#1} in the \OEISref.)}

\newcommand{\sub}[1]{\mathrm{Sub}(#1)}
\newcommand{\av}[1]{\mathrm{Av}(#1)}
\newcommand{\minbl}[3]{\mathrm{mb}(#1;#2,#3)}
\newcommand{\ext}[4]{\lrud{#1}_#2(#3,#4)}
\newcommand{\lgreedy}[2]{#1^{#2}_{\lambda}}
\newcommand{\rect}{\mathrm{rect}}
\newcommand{\pin}[1]{#1}

\begin{abstract}
A permutation class which is closed under pattern involvement may be
described in terms of its basis. The wreath product construction
$X\wr Y$ of two permutation classes $X$ and $Y$ is also closed, and
we investigate classes $Y$ with the property that, for any finitely
based class $X$, the wreath product $X\wr Y$ is also finitely based.
\end{abstract}



\section{Introduction and Statement of Theorem}

Two finite sequences of the same length, $\alpha=a_1a_2\cdots a_n$ 
and $\beta=b_1b_2\cdots b_n$, are said to be \emph{order isomorphic} 
if, for all $i,j$, we have $a_i < a_j$ if and only if $b_i < b_j$. 
Viewing permutations of length $n$ as orderings on the numbers 
$1,2,\ldots,n$, every sequence of $n$ distinct symbols is order 
isomorphic to a unique permutation. A permutation $\sigma$ is 
said to be \emph{involved} in the permutation $\pi$ (denoted 
$\sigma\leq\pi$) if there is a subsequence (or \emph{pattern}) 
of $\pi$ order isomorphic to $\sigma$\symbolfootnote[2]{For a sequence 
$\alpha$ (not necessarily a permutation) and set of permutations $Y$, with a slight
abuse of notation, we will sometimes write statements like
``$\alpha\in Y$'', meaning ``the permutation order isomorphic to
$\alpha$ lies in $Y$.''}. For example, $1324\leq 6351427$
because of the subsequence $3547$. A book introducing the study of these 
permutation patterns has been written by B\'{o}na~\cite{bona:combinatorics-o:}.

This involvement order forms a partial order on the set of all finite
permutations; sets of permutations which are closed downwards under
this order are called \emph{permutation classes}. These classes are specified primarily in one of three ways:
\begin{itemize}
\item {\bf Pattern avoidance.} A permutation class $X$ can be regarded as a set 
of permutations which \emph{avoid} certain patterns. The set $B$ of minimal permutations 
not in $X$ forms an antichain, and is known as the \emph{basis} of $X$. We write $X=\av{B}$ 
to mean the class $X=\{\pi \mid\beta\not\leq\pi \textrm{ for all } \beta\in B \}$. Antichains 
(and hence bases) need not be finite -- see, for example, Atkinson, Murphy and Ru\v{s}kuc~%
\cite{atkinson:partially-well-:}, Murphy~\cite{murphy:restricted-perm:} and Murphy and Vatter
\cite{murphy:profile-classes:}.

\item {\bf Permuting machines.} Permutation classes arise naturally as a 
result of machines which permute an input stream of symbols. The first 
such class to appear was the set of stack-sortable permutations, 
presented by Knuth~\cite{knuth:the-art-of-comp:}. 

\item {\bf Constructions.} New permutation classes can be formed using 
constructions involving one or more old classes. Atkinson~\cite{atkinson:restricted-perm:} 
gives the first study of these, and some further constructions 
can be found in Atkinson and Stitt~\cite{atkinson:restricted-perm:a} and 
Murphy~\cite{murphy:restricted-perm:}.
\end{itemize} 

In all but the first of these, a natural question to ask is if the class 
is finitely based. In the case of permuting machines -- more specifically, 
stack sorting -- B\'{o}na's survey~\cite{bona:a-survey-of-sta:} reviews several 
answers to this question. In the case of constructions, there are many
with only partial answers. Here, we will consider the question of basis 
for the wreath product, a construction which is intrinsically connected 
to simple permutations and the substitution decomposition 
-- see Albert and Atkinson \cite{albert:simple-permutat:} and Brignall, 
Huczynska and Vatter \cite{brignall:simple-permutat:}. A special case 
of the wreath product -- the ``profile classes'' of \cite{atkinson:restricted-perm:} 
-- was also used to give alternative proofs of the enumeration results 
in West \cite{west:generating-tree:}. 

Given a permutation $\pi\in S_n$ 
and nonempty permutations $\alpha_1,\alpha_2,\ldots,\alpha_n$, 
the \emph{inflation of $\pi$ by $\alpha_1,\alpha_2,\ldots,\alpha_n$} 
is the permutation obtained by replacing each point $\pi(i)$ by an interval 
order isomorphic to $\alpha_i$, and is denoted 
$\pi[\alpha_1,\alpha_2,\ldots,\alpha_n]$. For example, 
$132[21,2413,321] = 217968543$. Conversely, a
\emph{deflation} of $\pi$ is any permutation
$\sigma$ arising from a decomposition $\pi =
\sigma[\pi_1,\pi_2,\ldots,\pi_m]$.

The \emph{wreath product} of two sets of permutations $X$ and $Y$ (not
necessarily permutation classes) is the set $X\wr Y$ of
all permutations which can be expressed as an inflation of a
permutation in $X$ by permutations in $Y$, i.e.
the set of permutations of the form $\pi[\alpha_1,\alpha_2,\ldots,\alpha_n]$
with $\pi\in X$ and $\alpha_1,\alpha_2,\ldots,\alpha_n\in Y$. 
It is easy to check that the wreath product of two permutation classes is
again a permutation class, but in only a few cases is the question of finite basis answerable.  
It is proved in \cite{atkinson:restricted-perm:a} that for any finitely based class $X$, 
the wreath product $X\wr \av{21}$ is also finitely based, and 
that $\av{21}\wr\av{321654}$ is not finitely based. 
Our primary aim here is to establish the following general theorem:

\begin{thm}\label{thm-intro-wfbp}For any finitely based class $Y$ not 
admitting arbitrarily long pin sequences, the wreath product $X\wr Y$ 
is finitely based for all finitely based classes $X$.
\end{thm}

The approach is constructive; first we introduce $Y$-profiles, which
give us the ability to decompose permutations arising in wreath
products into components belonging to the two original classes. For
a permutation not arising in such a wreath product, we prove the
existence of a subsequence order isomorphic to a basis element of the 
class $X$. Moreover, there is a basis
element of $Y$ lying within the ``minimal block'' defined by any two 
points of this subsequence. It is then a matter of using these considerations to
show that, when the class $Y$ admits only finite pin sequences, the
minimal elements not in the wreath product have bounded size.

Our secondary aim, arising as a result of the above considerations, 
is to exhibit a number of classes of the form $Y=\av{\alpha}$ for $|\alpha|\leq 3$, 
or $Y=\av{\alpha,\beta}$ with $|\alpha|\leq 4$, $|\beta|\leq 4$ which do not satisfy 
Theorem~\ref{thm-intro-wfbp}, and to demonstrate how an infinitely based wreath 
product $X\wr Y$ can be found in each case.

\section{Simplicity and Substitution Decomposition}

As mentioned earlier, the wreath product is closely related to simple 
permutations and the substitution decomposition, both of which we will need, 
so here we review these concepts. Often we are going to view permutations 
as points in a plane; the \emph{plot} of a permutation $\pi$ is the set of coordinates
$\{(i,\pi(i))\}$ in the plane. This viewpoint will provide
invaluable insight into many of the structural considerations 
discussed later on.

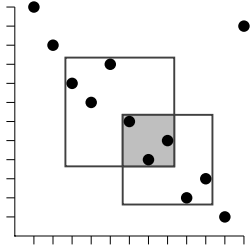
\begin{figure}
\begin{center}
\psset{xunit=0.01in, yunit=0.01in} \psset{linewidth=0.005in}
\begin{pspicture}(0,0)(120,120)
\psaxes[dy=10,dx=10,tickstyle=bottom,showorigin=false,labels=none](0,0)(120,120)
\psframe[linecolor=darkgray,fillcolor=lightgray,linewidth=0.01in](26,36)(84,94)
\psframe[linecolor=darkgray,linewidth=0.01in](56,16)(104,64)
\psframe[linecolor=darkgray,fillstyle=solid,fillcolor=lightgray,linewidth=0.01in](56,36)(84,64)
\pscircle*(10,120){0.03in} 
\pscircle*(20,100){0.03in} 
\pscircle*(30,80){0.03in}
\pscircle*(40,70){0.03in} 
\pscircle*(50,90){0.03in}
\pscircle*(60,60){0.03in} 
\pscircle*(70,40){0.03in}
\pscircle*(80,50){0.03in} 
\pscircle*(90,20){0.03in}
\pscircle*(100,30){0.03in}
\pscircle*(110,10){0.03in}
\pscircle*(120,110){0.03in}
\end{pspicture}
\end{center}
\caption{Two intervals and their intersection.}
\label{fig-intervals}
\end{figure}

An \emph{interval} or \emph{block} of a permutation $\pi$ is a
segment $\pi(i)\pi(i+1)\cdots \pi(i+j)$ in which the set of values forms 
an interval of natural numbers. In the plot of a permutation, intervals can 
be seen as a set of points enclosed in an axis-parallel rectangle, 
with no points lying in the regions above, below, to the left 
or to the right. It is worth noting that the intersection 
of two intervals is itself an interval, an observation clearly 
seen in Figure~\ref{fig-intervals}.

The permutation $\pi$ is \emph{simple} if its only intervals are 
singletons, or the whole of $\pi$. Note that simple permutations 
have only trivial deflations, and are the only permutations with 
this property. As such, they can be regarded as the building blocks 
of permutation classes. Every permutation can be written as the 
inflation of a unique simple permutation, and this decomposition 
is known as the \emph{substitution decomposition}. We shall refer 
to the unique simple permutation in this decomposition as the 
\emph{skeleton}. If the skeleton has length at least $4$, then the 
whole decomposition is unique:

\begin{prop}If $\pi$ has a substitution decomposition
$\sigma[\pi_1,\pi_2,\ldots,\pi_m]$ with $m\geq 4$, then every
$\pi_i$ is determined uniquely.\end{prop}

When $m=2$, we may write $\pi = 12[\pi_1,\pi_2]$, in which case
$\pi$ is \emph{sum decomposable}, or $\pi = 21[\pi_1,\pi_2]$, in
which case $\pi$ is \emph{skew decomposable}, and in both cases the
choice of $\pi_1,\pi_2$ is not necessarily unique. A permutation
that is not sum (respectively, skew) decomposable is \emph{sum}
(resp.~\emph{skew}) \emph{indecomposable}.

\section{\boldmath $Y$-Profiles}
We need to be able to know when a given permutation lies in 
the wreath product of two permutation classes. This could be done 
by inspecting all possible decompositions and checking for 
membership of the orignal classes, but this is liable to be 
computationally intensive. Instead, we would prefer only to 
check a single decomposition, from which membership or otherwise 
of the wreath product is immediately obvious.

The \emph{profile} of a permutation $\pi$ is the unique permutation
obtained by contracting every maximal consecutive increasing
sequence in $\pi$ into a single point~%
\cite{atkinson:restricted-perm:}. For example, the profile of
$3415672$ is $3142$ because of the segments 34, 1, 567 and 2.

The notion of a ``$Y$-profile'' connects this idea with the
definition of the substitution decomposition $\pi =
\sigma[\pi_1,\ldots,\pi_m]$ of $\pi$. We want the $Y$-profile of
$\pi$ to be the shortest possible deflation of $\pi$, given we may
only deflate by elements from the class $Y$. However, this is not
clearly well-defined, so before we can proceed, we must first
introduce $Y$-deflations.

Formally, let $Y$ be a permutation class, and $\pi$ any permutation.
Then a $Y$\emph{-deflation} of $\pi$ is a permutation $\pi'$ for
which $\pi$ can be expressed as
$\pi'[\alpha_1,\alpha_2,\ldots,\alpha_k]$ with
$\alpha_1,\alpha_2,\ldots,\alpha_k \in Y$. For an arbitrary permutation 
$\pi$, there are many different
$Y$-deflations. However, the shortest one is unique, and it is this
one that gives rise to the $Y$-profile.

\begin{lem}\label{yprofileunique}
For every closed class $Y$ and permutation $\pi$, the shortest
$Y$-deflation of $\pi$ is unique.
\end{lem}

\begin{proof} We proceed by induction on $n=|\pi|$. The case $n=1$ is
trivial, so now suppose $n>1$. Fix a shortest $Y$-deflation of the
permutation $\pi$, and label this permutation $\pi^Y$. If $\pi\in Y$
then $\pi^Y=1$ is unique, so we will assume $\pi\notin Y$.

Let $\sigma$, of length $m\geq 2$, be the skeleton of $\pi$, and first
consider the case where $m\geq 4$, whereby we have the unique
substitution decomposition $\pi = \sigma[\pi_1,\pi_2,\ldots,\pi_m]$.
By the inductive hypothesis, the shortest $Y$-deflations of
$\pi_1,\pi_2,\ldots,\pi_m$ are unique, and we will label them
$\pi_1^Y,\pi_2^Y,\ldots,\pi_m^Y$. We claim that
$\pi^Y=\sigma[\pi_1^Y,\pi_2^Y,\ldots,\pi_m^Y]$. Consider any other
$Y$-deflation of $\pi$, $\pi =
\pi'[\alpha_1,\alpha_2,\ldots,\alpha_k]$. Since $\pi\notin Y$,
$\pi'$ cannot be trivial, and so $\sigma\leq\pi'$, and indeed
$\sigma$ is the skeleton of $\pi'$, giving a unique deflation $\pi'
= \sigma[\pi'_1,\ldots,\pi'_m]$. Moreover, $\pi'_i$ is a
$Y$-deflation of $\pi_i$ for all $i$. Since $\pi_i^Y$ is the unique
shortest $Y$-deflation, we must have $\pi_i^Y\leq\pi'_i$, which
implies $\pi^Y\leq\pi'$.

When $m = 2$, more care is required. In this case $\pi$ is either
sum or skew decomposable, and without loss of generality we may
assume the former. Write $\pi = 12\cdots
t[\pi_1,\pi_2,\ldots,\pi_t]$ where each $\pi_i$ is sum
indecomposable. If every $\pi_i\in Y$, then any shortest
$Y$-deflation of $\pi$ will be an increasing permutation of length
at most $t$, and as there is only one increasing permutation of each
length, $\pi^Y$ will be unique. So now suppose that there exists at
least one $i$ such that $\pi_i\notin Y$, so that $|\pi_i^Y|\geq 2$.
Since $\pi_i$ is sum indecomposable, $\pi_i^Y$ is also sum
indecomposable. We claim the shortest $Y$-deflation of $\pi$ will be
$$\pi^Y = \left(\pi_1\oplus\cdots\oplus\pi_{i-1}\right)^Y \oplus
\pi_i^Y \oplus \left(\pi_{i+1}\oplus\cdots\oplus\pi_t\right)^Y.$$
Any other $Y$-deflation will also have to be written as a direct sum
of three permutations in this way, and by induction each of these
will involve the respective shortest $Y$-deflation.
\end{proof}

Thus, for any class $Y$ and permutation $\pi$, the
$Y$\emph{-profile} of $\pi$ is the unique shortest $Y$-deflation of
$\pi$, and is denoted $\pi^Y$. Note that setting $Y=\av{21}$, the
set of increasing permutations, returns the original definition of
the profile, but if we set $Y=S$, the set of all permutations, we do
not get the substitution decomposition back, as $\pi^S = 1$ for any
permutation. However, an easy consequence of the above proof is that
if $\pi\notin Y$, and $\sigma$ is the skeleton of $\pi$, then
$\sigma\leq\pi^Y$.

As mentioned at the beginning of this section, our aim with $Y$-profiles 
is to be able to to move from the permutations of the wreath product 
$X\wr Y$ down to the permutations in the two classes $X$ and $Y$ in a single step.
Thus although initially we may know very little about the structure of a
permutation in the basis of $X\wr Y$, by taking its $Y$-profile we
should be left with a permutation involving a (known) basis element
of $X$. Conversely, we want to be able to construct basis elements
of $X\wr Y$ given only the bases of $X$ and $Y$. These ideas are
encapsulated in the following theorem.

\begin{thm}\label{wreathprofilethm}
  Let $X$ and $Y$ be two arbitrary permutation classes. Then $\pi\in X\wr
  Y$ if and only if $\pi^{Y}\in X$.
\end{thm}

\begin{proof} One direction is immediate. For the converse, since $\pi\in
X\wr Y$, there exists $\pi' \in X$ which is a deflation of $\pi$ by
permutations in $Y$. The proof of Lemma~\ref{yprofileunique} then
tells us that $\pi^Y \leq \pi'$, completing the proof.\end{proof}

Any expression of the form $\pi=\pi^Y[\alpha_1,\ldots,\alpha_k]$ is
called a $Y$\emph{-profile decomposition} of $\pi$, and the blocks
$\alpha_i$ are called the \emph{$Y$-profile blocks}. These blocks
are not typically uniquely defined. For example, the
$\av{123}$-profile of $234615$ is $23514$, but it can be decomposed
either as $23514[12,1,1,1,1]$ or $23514[1,12,1,1,1]$. Thus it will
be useful to fix a particular $Y$-profile decomposition, especially
as later we are going to need to know about the structure of each
of the $Y$-profile blocks.

The \emph{left-greedy $Y$-profile} of $\pi$ is the decomposition
$\pi=\lgreedy{\pi}{Y}[\lambda_1,\lambda_2,\ldots,\lambda_\ell]$ with
$\lambda_i\in Y$ for all $i$, in which $\lambda_1$ is first chosen
maximally, then $\lambda_2$, and so on. Each $\lambda_i$ is called a
\emph{left-greedy $Y$-profile block} of $\pi$. This yields the
usual, unique, $Y$-profile:

\begin{lem}For any class $Y$ and permutation $\pi$, $\pi^Y = \lgreedy{\pi}{Y}$.\end{lem}

\begin{proof}
Again, we use induction on $n=|\pi|$. The base case $n=1$ is
trivial, so now suppose $n>1$. Assume further that $\pi\notin Y$, as
otherwise $\pi^Y=\lgreedy{\pi}{Y}=1$ follows immediately. Let
$\pi=\lgreedy{\pi}{Y}[\lambda_1,\lambda_2,\ldots,\lambda_\ell]$ be
the left-greedy $Y$-profile of $\pi$, let
$\pi^Y[\alpha_1,\alpha_2,\ldots,\alpha_k]$ be any other $Y$-profile
decomposition of $\pi$, and let $\sigma[\pi_1,\pi_2,\ldots,\pi_m]$
be the substitution decomposition.

Consider first the case where $m=|\sigma|\geq 4$. By the proof of
Lemma~\ref{yprofileunique}, we have
$\pi^Y=\sigma[\pi_1^Y,\pi_2^Y,\ldots,\pi_m^Y]$. A similar argument
shows that
$\lgreedy{\pi}{Y}=\sigma[\lgreedy{(\pi_1)}{Y},\lgreedy{(\pi_2)}{Y},\ldots,\lgreedy{(\pi_m)}{Y}]$,
and by induction $\pi_i^Y = \lgreedy{(\pi_i)}{Y}$ for all $i$,
giving the required result.

When $m=2$, $\pi$ is either sum or skew decomposable, and we may
assume the former. Write $\pi = 12\cdots
t[\pi_1,\pi_2,\ldots,\pi_t]$ where each $\pi_i$ is sum
indecomposable. In the case where every $\pi_i\in Y$, both $\pi^Y$
and $\lgreedy{\pi}{Y}$ will be increasing permutations with $k\leq
\ell\leq t$. When using the left-greedy $Y$-profile decomposition,
the block $\lambda_1$ was chosen maximally, and so $\alpha_1\leq
\lambda_1$. Then the block $\lambda_2$ was taken maximally, so the
$Y$-profile block $\alpha_2$ cannot extend further right than the
end of $\lambda_2$, hence $\alpha_2\leq
\lambda_1\oplus\lambda_2$. Continuing in this manner, we see that,
for all $i$, $\alpha_i\leq
\lambda_1\oplus\lambda_2\oplus\cdots\oplus\lambda_i$, and in
particular
$\alpha_k\leq\lambda_1\oplus\lambda_2\oplus\cdots\oplus\lambda_k$.
But we must have $k\leq\ell$, and so $k=\ell$. The remaining case is
where at least one $\pi_i\notin Y$. Pick $i$ to be minimal with this
property, and then by the proof of Lemma~\ref{yprofileunique},the
$Y$-profile breaks into three pieces,
$$\pi^Y = \left(\pi_1\oplus\cdots\oplus\pi_{i-1}\right)^Y \oplus \pi_i^Y \oplus \left(\pi_{i+1}\oplus\cdots\oplus\pi_t\right)^Y.$$
A similar argument holds for the left-greedy $Y$-profile, and then
by induction each of the three pieces in the left-greedy $Y$-profile
is equal to the corresponding piece in the $Y$-profile.\end{proof}

There is, of course, nothing special about the left-greedy
$Y$-profile; it can be seen that any algorithm to compute a
$Y$-profile-like decomposition in which at each stage the blocks are
chosen maximally will yield a $Y$-profile deflation. For our
purposes, however, when required we will always use the left-greedy
algorithm.

\section{The Minimal Block}
\begin{figure}
\begin{center}
\psset{xunit=0.01in, yunit=0.01in} \psset{linewidth=0.005in}
\begin{pspicture}(0,0)(180,168)
\psaxes[dy=15,dx=15,tickstyle=bottom,showorigin=false,labels=none](0,0)(165,165)
\psframe[linecolor=darkgray,fillstyle=solid,fillcolor=lightgray,linewidth=0.01in](23,38)(97,112)
\pscircle*(15,30){0.04in} \pscircle(30,45){0.04in}
\pscircle(45,90){0.04in} \pscircle*(60,105){0.04in}
\pscircle*(75,60){0.04in} \pscircle*(90,75){0.04in}
\pscircle*(105,135){0.04in} \pscircle*(120,120){0.04in}
\pscircle*(135,15){0.04in}
\end{pspicture}
\end{center}
\caption{The minimal block $\minbl{\pi}{2}{3}$ in $\pi=236745981$.}
\label{fig-min-block}
\end{figure}
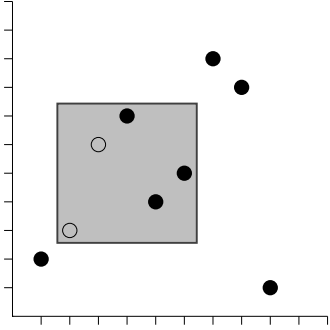

The primary aim of this section is to be able to tell if any
two points in a permutation belong to the same left-greedy
$Y$-profile block, and also a partial converse: given the
$Y$-profile deflation, what can we say about the points ``between''
two specified points? To this end, we define a new concept as
follows. Let $\pi$ be any permutation of length $n$. For all $1\leq
i < j \leq n$, the \emph{minimal block of $\pi$ that contains
$\pi(i)$ and $\pi(j)$}, denoted $\minbl{\pi}{i}{j}$, is the set of
points of $\pi$ which forms the shortest interval involving both
$\pi(i)$ and $\pi(j)$. In other words, there exists $k\leq i$ and
$\ell\geq j - k$ such that $\minbl{\pi}{i}{j}=\pi(k)\cdots
\pi(k+\ell)$ forms an interval but no subsegment of this contains
both $\pi(i)$ and $\pi(j)$ and forms an interval. For example, if
$\pi = 236745981$, then the minimal block on $\pi(2)=3$ and
$\pi(3)=6$ is $\minbl{\pi}{2}{3} = 36745$ (See Figure~%
\ref{fig-min-block}).

It follows from the observation that the intersection of two
intervals itself forms an interval that the minimal block is always
uniquely defined. Before we can proceed to the main result, we make
one further observation.

\begin{lem}\label{minimalblockinvolvementlemma}Let $\pi$ be any permutation and
  let $i\neq j$ be any pair of positions in $\pi$. Then if
  $k,\ell \in\minbl{\pi}{i}{j}$ with $k\neq\ell$ we have
  $$\minbl{\pi}{k}{l} \subseteq \minbl{\pi}{i}{j}.$$
  Moreover, if both $i$ and $j$ separate $k$ from $\ell$ by position, then $\minbl{\pi}{k}{\ell} = \minbl{\pi}{i}{j}$.
\end{lem}

\begin{proof} That $\minbl{\pi}{k}{\ell}$ is contained in
$\minbl{\pi}{i}{j}$ is obvious. Now suppose $i$ and $j$ separate $k$
from $\ell$ by position, i.e. $k\leq i<j \leq \ell$. Then
$\minbl{\pi}{k}{\ell}$ is an interval of $\pi$ involving both
$\pi(i)$ and $\pi(j)$. As $\minbl{\pi}{i}{j}$ is minimal with this
property, we have $\minbl{\pi}{i}{j}\subseteq\minbl{\pi}{k}{\ell}$
and so $\minbl{\pi}{i}{j}=\minbl{\pi}{k}{\ell}$.\end{proof}

We are now ready to prove our main technical result of this section.

\begin{lem}\label{minimalblocklem}
Let $Y$ be a permutation class, and let $\pi\in S_n$ be any
permutation. Then for any pair $i,j$ with $1\leq i<j\leq n$:
\begin{enumerate}
\item[(i)] If the permutation order isomorphic to $\minbl{\pi}{i}{j}$ does not lie in $Y$,
 then $\pi(i)$ and $\pi(j)$ lie in different $Y$-profile blocks.
\item[(ii)] Conversely, if $\pi(a_i)$ and $\pi(a_j)$ are the first symbols
  of two distinct left greedy $Y$-profile blocks $\alpha_i$ and $\alpha_j$ respectively,
  then the permutation order isomorphic to $\minbl{\pi}{i}{j}$ does not lie in $Y$.
\end{enumerate}
\end{lem}

\begin{proof} (i) By minimality and uniqueness of the minimal
  block, every block in $\pi$ containing both $\pi(i)$ and $\pi(j)$
  must contain the minimal block $\minbl{\pi}{i}{j}$. Hence every such block
  does not lie in $Y$, so cannot be a $Y$-profile block.

(ii) Write $\pi = \pi^{Y}[\alpha_1,\alpha_2,\ldots,\alpha_k]$, and
let the sequence $\pi(a_1),\pi(a_2),\ldots,\pi(a_k)$ represent the
leading points in $\pi$ of the left-greedy $Y$-profile blocks
$\alpha_1,\alpha_2,\ldots,\alpha_k$. Let $\alpha_i$ and $\alpha_j$,
$i<j$, be a pair of $Y$-profile blocks. We prove the statement by
induction on $i$.

 When $i=1$, the block $\alpha_1$ was picked maximally subject to
 $\alpha_1\in Y$. For any $j>1$, the minimal block $\minbl{\pi}{a_1}{a_j}$
 strictly contains $\alpha_1$ and then the maximality of $\alpha_1$ is
 contradicted unless $\minbl{\pi}{a_1}{a_j}\notin Y$.

Suppose now that $i>1$, and that $\minbl{\pi}{a_\ell}{a_j}\notin Y$
for any $\ell<i$ and $j>\ell$. The $Y$-profile block $\alpha_i$ was
picked maximally to avoid basis elements of $Y$, subject to starting
at symbol $\pi(a_i)$. Consider, for some $j>i$, the minimal block
$\minbl{\pi}{a_i}{a_j}$, necessarily containing all of $\alpha_i$.
If the leftmost point of $\minbl{\pi}{a_i}{a_j}$ is $\pi(a_i)$, then
since $\alpha_i$ is the maximal block lying in $Y$
  which starts at $\pi(a_i)$, we must have $\minbl{\pi}{a_i}{a_j} \notin Y$.
  So now suppose that $\minbl{\pi}{a_i}{a_j}$ contains at least one symbol
  $\pi(h)$ from $\pi$ with $h<a_i$. Let the $Y$-profile block containing $\pi(h)$
  be $\alpha_\ell$; we claim that $\alpha_\ell$ is completely contained in $\minbl{\pi}{a_i}{a_j}$.
  If not, then part of $\alpha_\ell$ lies outside $\minbl{\pi}{a_i}{a_j}$
  in both position and value, and so the part lying inside $\minbl{\pi}{a_i}{a_j}$
  itself forms an interval in either the top-left or bottom-left corner of
  the minimal block, but yet it contains neither $\pi(a_i)$ nor $\pi(a_j)$, 
  contradicting the minimality of $\minbl{\pi}{a_i}{a_j}$. In particular, 
  the first symbol $\pi(a_\ell)$ of $\alpha_\ell$
  is in $\minbl{\pi}{a_i}{a_j}$, and by Lemma~\ref{minimalblockinvolvementlemma},
  we have $\minbl{\pi}{a_\ell}{a_j} = \minbl{\pi}{a_i}{a_j}$. By the inductive
  hypothesis $\minbl{\pi}{a_\ell}{a_j}\notin Y$, and so $\minbl{\pi}{a_i}{a_j}\notin Y$.\end{proof}

Using this result, we now know when two points of a permutation will
lie in the same $Y$-profile block, and, more importantly for what
follows, we know that a basis element of $Y$ exists in the minimal
block of the first symbols of any two $Y$-profile blocks. What we do
not yet know is how to find it; given such a minimal block, we need
a method to search through the block systematically and locate the
points that form this basis element within a bounded number of
steps. This is the subject of the next section.

\section{Pin Sequences and the Wreath Product}

Pin sequences were introduced by Brignall, Huczynska and Vatter~%
\cite{brignall:simple-permutat:a} in the study of simple permutations.
The idea there is that, since simple permutations have no
non-trivial intervals, if we begin with any two points we can use
pin sequences to get to any chosen edge of our simple permutation.
Here we are not working solely with simple permutations, and we cannot
expect the same result to hold in the general case. However, we can
obtain the same result for the minimal block. We begin by reviewing
some terminology from \cite{brignall:simple-permutat:a}, and to do this 
it is best to revert to viewing permutations as plots in the plane. 

\begin{figure}
\begin{center}
\psset{xunit=0.01in, yunit=0.01in}
\psset{linewidth=0.005in}
\begin{pspicture}(0,0)(180,168)
\psaxes[dy=15,dx=15,tickstyle=bottom,showorigin=false,labels=none](0,0)(165,165)
\psframe[linecolor=darkgray,fillstyle=solid,fillcolor=lightgray,linewidth=0.02in](57,27)(168,138)
\psline(75,75)(120,75)
\psline(105,60)(105,135)
\psline(90,90)(135,90)
\psline(90,120)(165,120)
\psline(150,135)(150,30)
\psline(15,45)(165,45)
\pscircle*(15,45){0.04in}
\pscircle*(30,150){0.04in}
\pscircle*(45,15){0.04in}
\pscircle*(60,105){0.04in}
\pscircle*(75,165){0.04in}
\pscircle*(90,60){0.04in}
\pscircle*(105,135){0.04in}
\pscircle*(120,75){0.04in}
\pscircle*(135,90){0.04in}
\pscircle*(150,30){0.04in}
\pscircle*(165,120){0.04in}
\rput[c](75,105){$p_1$}
\rput[c](75,60){$p_2$}
\rput[c](135,67){$p_3$}
\rput[c](105,148){$p_4$}
\rput[c](135,103){$p_5$}
\rput[c](180,120){$p_6$}
\rput[c](150,17){$p_7$}
\rput[c](15,58){$p_8$}
\end{pspicture}
\end{center}
\caption{A pin sequence.}
\label{fig-pin-sequence}
\end{figure}
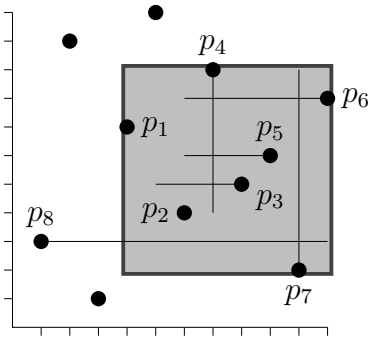

For points $\pin{p}_1,\pin{p}_2,\ldots,\pin{p}_m$ in the plane, let
$\rect(\pin{p}_1,\pin{p}_2,\ldots,\pin{p}_m)$ be the smallest
axis-parallel rectangle containing them. Note that this is different
to the minimal block, as we do not require that
$\rect(\pin{p}_1,\pin{p}_2,\ldots,\pin{p}_m)$ be an interval.

Let $\pi$ be a permutation. A \emph{pin sequence} is a sequence of points $\pin{p}_1,\pin{p}_2,\ldots$ of $\pi$ which for $i\geq 3$ obey, when plotted in a plane, the following two conditions.
\begin{itemize}
\item $\pin{p}_i\notin\rect(\pin{p}_1,\pin{p}_2,\ldots,\pin{p}_{i-1})$,
\item $\pin{p}_i$ \emph{slices} $\rect(\pin{p}_1,\pin{p}_2,\ldots,\pin{p}_{i-1})$ either horizontally or vertically. That is $\pin{p}_i$ lies between two points of $\rect(\pin{p}_1,\pin{p}_2,\ldots,\pin{p}_{i-1})$ either by position or value.
\end{itemize}

For each pin $\pin{p}_i$, $i\geq 3$, we also specify a direction,
being left, right, up or down. For example, a \emph{left pin} is one
that lies between two point of
$\rect(\pin{p}_1,\pin{p}_2,\ldots,\pin{p}_{i-1})$ by value, but
whose position is smaller than any point of
$\rect(\pin{p}_1,\pin{p}_2,\ldots,\pin{p}_{i-1})$. In Figure
\ref{fig-pin-sequence}, $\pin{p}_3,\pin{p}_5$ and $\pin{p}_6$ are
right pins, $\pin{p}_4$ is an up pin, $\pin{p}_7$ a down pin and
$\pin{p}_8$ a left pin.

We create a \emph{proper pin sequence} by adjoining two further conditions:
\begin{itemize}
\item \emph{Maximality}: each pin must be taken maximally in its direction. For example, a proper left pin out of $\rect(\pin{p}_1,\pin{p}_2,\ldots,\pin{p}_{i-1})$ must be the left pin slicing $\rect(\pin{p}_1,\pin{p}_2,\ldots,\pin{p}_{i-1})$ with smallest position.
\item \emph{Separation}: in slicing $\rect(\pin{p}_1,\pin{p}_2,\ldots,\pin{p}_i)$, $\pin{p}_{i+1}$ must lie between $\pin{p}_{i}$ and
$\rect(\pin{p}_1,\pin{p}_2,\ldots,\pin{p}_{i-1})$ either by position or value.
\end{itemize}
For example, in Figure~\ref{fig-pin-sequence}, $\pin{p}_8$ is a proper left pin as it slices $\pin{p}_7$ from $\rect(\pin{p}_1,\pin{p}_2,\ldots,\pin{p}_6)$ and is maximal in its direction. Similarly, $\pin{p}_4$ and $\pin{p}_7$ are proper pins, but $\pin{p}_3, \pin{p}_5$ and $\pin{p}_6$ are not, as $\pin{p}_3$ does not obeying maximality, $\pin{p}_5$ does not separate $\pin{p}_4$ from $\rect(\pin{p}_1,\pin{p}_2,\pin{p}_3)$, and $\pin{p}_6$ does not separate $\pin{p}_5$ from $\rect(\pin{p}_1,\pin{p}_2,\pin{p}_3,\pin{p}_4)$.

In a proper pin sequence, the maximality and separation conditions force the pin $\pin{p}_{i+1}$ to have direction perpendicular to the direction of $\pin{p}_i$, so for example a left pin can only be followed by an up pin or a down pin.

If a pin sequence $\pin{p}_1,\pin{p}_2,\ldots,\pin{p}_m$ of $\pi$ is
such that $\rect(\pin{p}_1,\pin{p}_2,\ldots,\pin{p}_m)$ encloses all
of $\pi$, then we say that it is \emph{saturated}. When we restrict
to proper pin sequences this is likely to be impossible to acheive,
even in simple permutations. However a weaker condition does hold. A
pin sequence $\pin{p}_1,\pin{p}_2,\ldots,\pin{p}_m$ of $\pi$ is said
to be \emph{right-reaching} if
$\rect(\pin{p}_1,\pin{p}_2,\ldots,\pin{p}_m)$ encloses all of $\pi$:

\begin{prop}[Brignall, Huczynska, Vatter~\cite{brignall:simple-permutat:a}]\label{rightreaching}
From any pair of points in a simple permutation, there exists a
proper right-reaching pin sequence.
\end{prop}

Since we are not working solely with simple permutations, we need to modify this proposition. Instead, we want the same to hold within a minimal block, defined as usual by two points, which also form the first two points of our proper pin sequence. Here, right-reaching means that the last pin is the right-most point of the minimal block, rather than of the whole permutation. Hence:

\begin{lem}\label{minblrightreaching}Let $\pi\in S_n$ be any permutation, and let $1\leq i<j\leq n$. Then
there exists a proper pin sequence with starting points
$\pin{p}_1=(i,p_i)$ and $\pin{p}_2=(j,p_j)$ which is right-reaching
in $\minbl{\pi}{i}{j}$.\end{lem}

\begin{proof}In the minimal block $\minbl{\pi}{i}{j}$, there exists a saturated (non-proper)
pin sequence $\pin{p}_1,\pin{p}_2,\ldots$ starting from the pins
$\pin{p}_1=(i,\pi(i))$ and $\pin{p}_2=(j,\pi(j))$. If there were no
such sequence, then some corner of the minimal block, not including
either $\pi(i)$ or $\pi(j)$, would form an interval by itself,
contradicting the minimality of $\minbl{\pi}{i}{j}$. Moreover, we
may assume, by removing unnecessary pins and relabelling, that every
pin is maximal in its direction.

The proof then follows the proof in
\cite{brignall:simple-permutat:a} of Proposition~%
\ref{rightreaching}. Since the pin sequence is saturated, it
includes the rightmost point of $\pi$. Label this point
$\pin{p}_{i_1}$. Next, take the smallest $i_2<i_1$ such that
$\pin{p}_1,\pin{p}_2,\ldots,\pin{p}_{i_2},\pin{p}_{i_1}$ is a valid
pin sequence, and observe that $\pin{p}_{i_1}$ separates
$\pin{p}_{i_2}$ from
$\rect(\pin{p}_1,\pin{p}_2,\ldots,\pin{p}_{i_2-1})$, as
$\pin{p}_1,\pin{p}_2,\ldots,\pin{p}_{i_2-1},\pin{p}_{i_1}$ is not a
valid pin sequence. Continue in this manner, finding pins
$\pin{p}_{i_3},\pin{p}_{i_4},\ldots$ until we reach
$\pin{p}_{i_{m+1}} = \pin{p}_2$, and then
$\pin{p}_1,\pin{p}_2,\pin{p}_{i_m},\pin{p}_{i_{m-1}}\ldots,\pin{p}_{i_1}$
is a proper right-reaching pin sequence.\end{proof}

Proposition~\ref{rightreaching} is easily recovered from Lemma~%
\ref{minblrightreaching} by setting $\pi$ to be a simple
permutation, and observing that all minimal blocks in a simple permutation are the whole permutation.

We are now ready to prove our main result.

\begin{thm}\label{wfbpthm}Let $Y = \av{B}$ be a finitely based permutation class not admitting arbitrarily long pin sequences. Then $X\wr Y$ is finitely based for all finitely based classes $X=\av{D}$.
\end{thm}

\begin{proof} Let $b= \max_{\beta\in B}(|\beta|)$, $d = \max_{\delta\in D}(|\delta|)$, and $\pi$
be any permutation in the basis of $X\wr Y$. By Theorem~%
\ref{wreathprofilethm}, we have $\pi^Y\notin X$, and so there exists
some $\delta\in D$ such that $\delta \leq\pi^{Y}$. We will be
done if we can identify a bounded subsequence of $\pi$ order isomorphic to 
a permutation $\omega$, say, for which $\delta\leq\omega^Y$, as then
$\omega^Y\notin X$ implies $\omega\notin X\wr Y$, and hence
$\omega=\pi$.

First include in our subsequence of $\pi$ the set of points order 
isomorphic to $\delta$ with positions $d_1,d_2,\ldots,d_k$ ($k=|\delta|$), chosen 
so that each $\pi(d_i)$ is the leftmost point of a distinct left 
greedy $Y$-profile block, and the choice
of blocks is also leftmost. For every pair $d_i,d_{i+1}$, 
Lemma~\ref{minimalblocklem} tells us that the minimal block $\minbl{\pi}{d_i}{d_{i+1}}$ 
involves some $\beta\in B$, and we include one such occurrence of this $\beta$ in
our subsequence. Our aim now is to add a bounded number of points 
so that $\beta$ still lies in the minimal block of the 
permutation $\omega$ on the points corresponding to 
$\pi(d_i)$ and $\pi(d_{i+1})$, as then these two points 
are preserved distinctly in $\omega^Y$.
We do this by taking a proper right-reaching and a proper 
left-reaching pin sequence of $\minbl{\pi}{d_i}{d_{i+1}}$ 
(which exist by Lemma~\ref{minblrightreaching}), 
and including them in the subsequence. These pin sequences
are only guaranteed to be bounded when $Y$ does not 
admit arbitrarily long pin sequences, as then there 
exists a number $N$ so that every pin sequence of 
length $N+2$ involves some basis element of $Y$.

Thus $\omega^Y$ still involves a subsequence order isomorphic to
$\delta$, and $|\omega|\leq d + (d-1)(2(N-1)+b)$.
\end{proof}

Brignall, Ru\v{s}kuc and Vatter~\cite{brignall:simple-permutat:b}
proved that determining whether a finitely based class does not admit arbitrarily long pin
sequences is decidable, and therefore given any pattern class we
can tell whether Theorem~\ref{wfbpthm} applies.


\section{Infinitely Based Examples}

For a class $Y$ which admits infinite pin sequences, Theorem~%
\ref{wfbpthm} gives us no information on whether the basis of $X\wr
Y$ (here for a specified class $X$) is finite. However, the proof
does tell us what some of the basis elements look like, namely
permutations built around a basis element of $X$, and in the minimal 
block between each pair of these points, there is a basis element of $Y$.
Constructing arbitrarily long basis elements of this
type is then achieved by embedding arbitrarily long pin sequences in
the minimal blocks. For example, the class $\av{321}$ admits the
infinite pin sequence formed by alternating between up and right
pins, and so we have:

\begin{thm}\label{thm-25134-basis}$\av{25134}\wr\av{321}$ is not finitely based.\end{thm}

\begin{proof}
We exhibit an antichain generated by repeatedly taking up and right
pins lying in the basis of $\av{25134}\wr\av{321}$. The first few
elements of the antichain are
\begin{eqnarray*}\beta_1 &=& 2,5,1,3,7,6,4\\
\beta_2 &=& 2,5,1,3,7,4,9,8,6\\
\beta_k &=& 2,5,1,3,7,4 \mid 9,6,11,8,\ldots,2k+3,2k\mid
2k+5,2k+4,2k+2 \quad(k\geq 3).\end{eqnarray*}
\begin{figure}
\begin{center}
\psset{xunit=0.01in, yunit=0.01in} \psset{linewidth=0.005in}
\begin{pspicture}(0,0)(155,155)
\psaxes[dy=10,dx=10,tickstyle=bottom,showorigin=false,labels=none](0,0)(150,150)
\psccurve[linewidth=0.01in](145,135)(130,140)(125,155)(140,150)
\psccurve[linewidth=0.01in](3,20)(15,55)(65,45)(30,3)
\psline(140,110)(140,134)\psline(150,120)(100,120)
\psline(110,130)(110,90) \psline(120,100)(80,100)
\psline(90,110)(90,70) \psline(100,80)(60,80) \psline(70,50)(70,90)
\psline(80,60)(40,60) \psline(50,70)(50,20)
\pscircle*(150,120){0.03in} \pscircle*(140,140){0.03in}
\pscircle*(130,150){0.03in} \pscircle*(120,100){0.03in}
\pscircle*(110,130){0.03in} \pscircle*(100,80){0.03in}
\pscircle*(90,110){0.03in} \pscircle*(80,60){0.03in}
\pscircle*(70,90){0.03in} \pscircle*(60,40){0.03in}
\pscircle*(50,70){0.03in} \pscircle*(40,30){0.03in}
\pscircle*(30,10){0.03in} \pscircle*(20,50){0.03in}
\pscircle*(10,20){0.03in}
\end{pspicture}
\end{center}
\caption{The element $\beta_5$ in the basis of
$\av{25134}\wr\av{321}$.} \label{fig-321-antichain}
\end{figure}
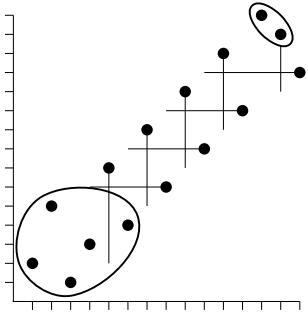
Here, as in \cite{atkinson:partially-well-:}, the $\mid$ symbol is
used only to clarify the structure of the permutation. See
Figure~\ref{fig-321-antichain} for an illustration of a typical
member of this antichain. We observe:
\begin{enumerate}
\item[(i)] The set $\{\beta_k\mid k\geq 1\}$ is an antichain.
\item[(ii)] The only occurrence of $321$ in each $\beta_k$ is $2k+5,2k+4,2k+2$.
\item[(iii)] The only occurrence of $25134$ in each $\beta_k$ is
  $2,5,1,3,\cdot,4$.
\item[(iv)] Each $\beta_k$ is neither sum nor skew decomposable.
\item[(v)] The $\av{321}$-profile of $\beta_k$ is
$2,5,1,3,7,4,\ldots,2k+3,2k,2k+4,2k+2$ (the only nontrivial
deflation occurs between $2k+5$ and $2k+4$). In particular,
$25134\prec\beta_k^{\av{321}}$ for all $k$, hence by Theorem~%
\ref{wreathprofilethm} $\beta_k\notin \av{25134}\wr\av{321}$.
\end{enumerate}
It only remains to show that $\beta_k$ is minimally not in
$\av{25134}\wr\av{321}$. Consider the effect of removing any symbol
$j$. If $j=2k+5, 2k+4$ or $2k+2$ then by (ii) this no longer
involves
  321 so $\beta_k - j \in \av{321} \subset \av{25134}\wr\av{321}$.
Similarly, if $j=2,5,1,3$ or $4$ then by (iii) $\beta_k-j$ no longer
involves
  $25134$ so $\beta_k - j \in \av{25134} \subset \av{25134}\wr\av{321}$.

For any other $j$, $\beta_k-j$ is sum decomposable. Under
  the $\av{321}$-profile, the first (lower) component deflates to a single point, and
  hence $(\beta_k-j)^{\av{321}}\in\av{25134}$. Thus $\beta_k - j \in
  \av{25134}\wr\av{321}$, completing the proof.\end{proof}

Note that in the above example, the class $X=\av{25134}$ was
specifically chosen so that the basis element $25134$ is not
contained in the repeated pin sequence used to build the antichain,
but it does lie in the class $Y$. This ensures that 25134 acts as an
``anchor'' at the base of the antichain, but yet the only instance
of the basis element 321 is in the upper ``anchor''.

As a result, for any class $Y$ which contains both the infinite pin
sequence formed by alternating between up and right pins, and the
permutation $25134$, the wreath product $\av{25134}\wr Y$ will
always contain an infinite antichain similar to the one above.
\begin{exa}
\begin{enumerate}
\item[(i)] The classes $Y=\av{321,2341}$ and $Y = \av{321,3412}$
both avoid the permutation $321$ and so the antichain in the proof
of Theorem~\ref{thm-25134-basis} lies in the basis of $\av{25134}\wr
Y$ in both cases.

\item[(ii)] All of the classes $Y=
\av{\alpha,\beta}$ with $(\alpha,\beta)$ being $(4321,4312)$,
$(4321,4231)$, $(4321,4213)$, $(4321,3412)$ and $(4321,3214)$ avoid
$4321$, and so the antichain with terms
\begin{eqnarray*}
\beta_1 &=& 2,5,1,3,8,7,6,4\\
\beta_2 &=& 2,5,1,3,7,4,10,9,8,6\\
\beta_k &=&2,5,1,3,7,4\mid 9,6,11,8,\ldots,2k+3,2k\mid
2k+6,2k+5,2k+4,2k+2 \quad(k\geq 3)\end{eqnarray*} lies in the basis
of $\av{25134}\wr Y$ in each case.

\item[(iii)] The classes $Y=\av{4312, 4231}$, $Y=\av{4312, 4213}$ and
$Y=\av{4312,3421}$ all avoid $4312$, so reversing the final two
points of each $\beta_k$ in case (ii) gives the required antichain.
\end{enumerate}
\end{exa}

\begin{exa} The two classes $Y=\av{4321,4123}$ and
$Y=\av{4312,4123}$ both admit the pin sequence formed by repeatedly
taking up and right pins, but do not contain the permutation
$25134$, because of the basis element $4123$. However, the class
$X=\av{25143}$ may be used instead. In the first case, the antichain
is (see Figure \ref{fig-4321-antichain} for an illustration):
\begin{figure}
\begin{center}
\psset{xunit=0.01in, yunit=0.01in} \psset{linewidth=0.005in}
\begin{pspicture}(0,0)(165,165)
\psaxes[dy=10,dx=10,tickstyle=bottom,showorigin=false,labels=none](0,0)(160,160)
\psccurve[linewidth=0.01in](125,165)(147,157)(155,135)(133,143)
\psccurve[linewidth=0.01in](3,20)(15,55)(65,35)(30,3)
\psline(150,110)(150,133)\psline(160,120)(100,120)
\psline(110,130)(110,90) \psline(120,100)(80,100)
\psline(90,110)(90,70) \psline(100,80)(60,80) \psline(70,50)(70,90)
\psline(80,60)(40,60) \psline(50,70)(50,20)
\pscircle*(160,120){0.03in} \pscircle*(150,140){0.03in}
\pscircle*(140,150){0.03in} \pscircle*(130,160){0.03in}
\pscircle*(120,100){0.03in} \pscircle*(110,130){0.03in}
\pscircle*(100,80){0.03in} \pscircle*(90,110){0.03in}
\pscircle*(80,60){0.03in} \pscircle*(70,90){0.03in}
\pscircle*(60,30){0.03in} \pscircle*(50,70){0.03in}
\pscircle*(40,40){0.03in} \pscircle*(30,10){0.03in}
\pscircle*(20,50){0.03in} \pscircle*(10,20){0.03in}
\end{pspicture}
\end{center}
\caption{The element $\beta_5$ in the basis of
$\av{25143}\wr\av{4321,4123}$.} \label{fig-4321-antichain}
\end{figure}
\begin{eqnarray*}
\beta_1 &=& 2,5,1,4,8,7,6,3\\
\beta_2 &=& 2,5,1,4,7,3,10,9,8,6\\
\beta_k &=& 2,5,1,4,7,3\mid 9,6,11,8,\ldots,2k+3,2k\mid
2k+6,2k+5,2k+4,2k+2 \quad(k\geq 3).\end{eqnarray*}
\end{exa}

All the examples so far have admitted the same ``up-right'' pin 
sequence. Another commonly found infinite pin sequence is 
formed by repeating the pattern left, down, right, up\symbolfootnote[1]{This
repeating pattern is the foundation for the ``Widdershins''
antichain of \cite{murphy:restricted-perm:}.}, and there are two classes
of the form $Y=\av{\alpha,\beta}$ with $|\alpha|=|\beta|=4$ which 
admit this sequence: $Y=\av{3412,2413}$ and $Y=\av{3412,2143}$. Each one 
must be handled separately.

\begin{exa}
\begin{enumerate}
\item[(i)]$Y=\av{3412,2413}$ may be paired with $X = \av{31542}$ to produce
  the antichain with terms
\begin{eqnarray*}
\beta_1&=&8,1,6,4,9,7,5,2,3\\
\beta_k&=&4k+4,1,4k+2,4,4k,6,\ldots 2k+6,2k\mid \\
&&2k+4,2k+2,2k+7,2k+5,2k+3\mid \\
&&2k+9,2k+1,\ldots,4k+5,5\mid 2,3 \quad(k\geq 2).
\end{eqnarray*}
See Figure~\ref{fig-widderschin-antichain} for an illustration. %
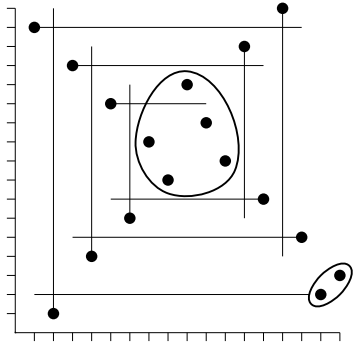
\begin{figure}
\begin{center}
\psset{xunit=0.01in, yunit=0.01in} \psset{linewidth=0.005in}
\begin{pspicture}(0,0)(175,175)
\psaxes[dy=10,dx=10,tickstyle=bottom,showorigin=false,labels=none](0,0)(170,170)
\psccurve[linewidth=0.01in](155,15)(170,20)(175,35)(160,30)
\psccurve[linewidth=0.01in](63,100)(80,73)(115,85)(90,137)
\psline(50,120)(100,120) \psline(60,130)(60,60)
\psline(50,70)(130,70) \psline(120,60)(120,150)
\psline(30,140)(130,140) \psline(40,150)(40,40)
\psline(30,50)(150,50) \psline(140,40)(140,170)
\psline(10,160)(150,160) \psline(20,170)(20,10)
\psline(10,20)(154,20)

\pscircle*(10,160){0.03in} \pscircle*(20,10){0.03in}
\pscircle*(30,140){0.03in} \pscircle*(40,40){0.03in}
\pscircle*(50,120){0.03in} \pscircle*(60,60){0.03in}

\pscircle*(70,100){0.03in} \pscircle*(80,80){0.03in}
\pscircle*(90,130){0.03in} \pscircle*(100,110){0.03in}
\pscircle*(110,90){0.03in}

\pscircle*(120,150){0.03in} \pscircle*(130,70){0.03in}
\pscircle*(140,170){0.03in} \pscircle*(150,50){0.03in}
\pscircle*(160,20){0.03in} \pscircle*(170,30){0.03in}
\end{pspicture}
\end{center}
\caption{The basis element $\beta_3$ in
$\av{31542}\wr\av{3412,2413}$.} \label{fig-widderschin-antichain}
\end{figure}%
Note that the occurrence of 3412 in any $\beta_k$ is not unique, but
every occurrence requires the final two symbols $2,3$ of $\beta_k$,
and so these points still behave in the same way as in previous
examples.

\item[(ii)]$Y=\av{3412,2143}$ may be paired with $X = \av{412563}$ to produce
  the antichain with terms:
\begin{eqnarray*}
\beta_1 &=& 10,1,8,4,6,9,11,7,5,2,3\\
\beta_k&=&4k+6,1,4k+4,4,4k+2,6,\ldots, 2k+8,2k \mid\\
&&2k+6,2k+2,2k+4,2k+7,2k+9,2k+5,2k+3\mid\\
&&2k+11,2k+1,\ldots,4k+7,5\mid 2,3 \quad(k\geq 2).
\end{eqnarray*}
\end{enumerate}\end{exa}

\section{Concluding Remarks}

The above examples suggest, to some extent, a general method for
finding infinite bases. However, these examples rely on just one
method for constructing antichains, and there is no reason why this
method should always work\symbolfootnote[3]{A somewhat different construction
was used by Atkinson and Stitt \cite{atkinson:restricted-perm:a} to
demonstrate an infinite antichain in the basis of
$\av{21}\wr\av{321654}$, relying on the sum decomposability of the
basis element $321654$.}. Moreover, within this construction,
finding a suitable class $X$ for a given class $Y$ is very specific
in each case.

In fact, it is unlikely that we can always find such a class $X$.
For example, the class of all subpermutations of the 
\emph{increasing oscillating sequence}, $416385\cdots$, is given by
$\av{321,2341,3412,4123}$
\cite{brignall:simple-permutat:b}, and admits the infinite proper pin
sequence alternating between an up pin and a right pin. However,
there are no other permutations in this class which can be used to 
anchor an infinite antichain based around this pin sequence,
so the method described hitherto does not work here. We therefore
pose the following question.

\begin{question}
Is there a finitely based class $X$ for which $X\wr
\av{321,2341,3412,4123}$ is not finitely based?
\end{question}

\paragraph{Acknowledgements.} The author wishes to thank Nik Ru\v{s}kuc 
and Vince Vatter for their invaluable comments.

\bibliographystyle{abbrv}
\bibliography{../refs}

\def\cprime{$'$}
\begin{thebibliography}{10}

\bibitem{albert:simple-permutat:}
M.~H. Albert and M.~D. Atkinson.
\newblock Simple permutations and pattern restricted permutations.
\newblock {\em Discrete Math.}, 300(1-3):1--15, 2005.

\bibitem{atkinson:restricted-perm:}
M.~D. Atkinson.
\newblock Restricted permutations.
\newblock {\em Discrete Math.}, 195(1-3):27--38, 1999.

\bibitem{atkinson:partially-well-:}
M.~D. Atkinson, M.~M. Murphy, and N.~Ru{\v{s}}kuc.
\newblock Partially well-ordered closed sets of permutations.
\newblock {\em Order}, 19(2):101--113, 2002.

\bibitem{atkinson:restricted-perm:a}
M.~D. Atkinson and T.~Stitt.
\newblock Restricted permutations and the wreath product.
\newblock {\em Discrete Math.}, 259(1-3):19--36, 2002.

\bibitem{bona:a-survey-of-sta:}
M.~B{\'o}na.
\newblock A survey of stack-sorting disciplines.
\newblock {\em Electron. J. Combin.}, 9(2):Article 1, 16 pp. (electronic),
  2003.

\bibitem{bona:combinatorics-o:}
M.~B{\'o}na.
\newblock {\em Combinatorics of permutations}.
\newblock Discrete Mathematics and its Applications (Boca Raton). Chapman \&
  Hall/CRC, Boca Raton, FL, 2004.
\newblock With a foreword by Richard Stanley.

\bibitem{brignall:simple-permutat:a}
R.~Brignall, S.~Huczynska, and V.~Vatter.
\newblock Decomposing simple permutations, with enumerative consequences.
\newblock \texttt{arXiv:math.CO/0606186}.

\bibitem{brignall:simple-permutat:}
R.~Brignall, S.~Huczynska, and V.~Vatter.
\newblock Simple permutations and algebraic generating functions.
\newblock \texttt{arXiv:math.CO/0608391}.

\bibitem{brignall:simple-permutat:b}
R.~Brignall, N.~Ru\v{s}kuc, and V.~Vatter.
\newblock Simple permutations: decidability and unavoidable substructures.
\newblock \texttt{arXiv:math.CO/0609211}.

\bibitem{knuth:the-art-of-comp:}
D.~E. Knuth.
\newblock {\em The art of computer programming. {V}ol. 1: {F}undamental
  algorithms}.
\newblock Addison-Wesley Publishing Co., Reading, Mass., 1969.

\bibitem{murphy:restricted-perm:}
M.~M. Murphy.
\newblock {\em Restricted permutations, antichains, atomic classes, and stack
  sorting}.
\newblock PhD thesis, Univ. of St Andrews, 2002.

\bibitem{murphy:profile-classes:}
M.~M. Murphy and V.~Vatter.
\newblock Profile classes and partial well-order for permutations.
\newblock {\em Electron. J. Combin.}, 9(2):Research paper 17, 30 pp.
  (electronic), 2003.

\bibitem{west:generating-tree:}
J.~West.
\newblock Generating trees and forbidden subsequences.
\newblock {\em Discrete Math.}, 157(1-3):363--374, 1996.

\end{thebibliography}

\end{document}